\documentclass[11pt]{amsart}
\usepackage{fullpage}
\usepackage{bm}
\usepackage{amsfonts,amssymb,mathrsfs,graphicx,color,xypic}
\usepackage[latin1]{inputenc}
\usepackage[english]{babel}
 \xyoption{all}

 \newcommand{\F}{\mathbb{F}}

 \theoremstyle{plain}
 \newtheorem{thm}{Theorem}[section]
 \newtheorem{defi}[thm]{Definition}

 \newtheorem{prop}[thm]{Proposition}
 \newtheorem{lem}[thm]{Lemma}
 \newtheorem{coro}[thm]{Corollary}

 \theoremstyle{remark}
 \newtheorem{rem}[thm]{Remark}

 \newenvironment{preuve}{{\bf Proof.}}{\hfill$\square$}

\date{March 2009}

\begin{document}

 %%%%% ------------- fill in your data below this line
 %%%%%-------------------
 %%%%%    The following lines \Title ... \EndAddress must ALL be present
 %%%%%    and in the given order.

\title{Completely symmetric configurations for $\sigma$-games on grid graphs}
\author{Mathieu Florence}
\address{Universit\'e Paris 6, Equipe de Topologie et G\'eom\'etrie Alg\'ebriques, bureau 9D07, 175 rue du Chevaleret, 75013 Paris.}
\email{mathieu.florence@gmail.com}
\author{Fr\'ed\'eric Meunier}
\address{Universit\'e Paris Est, LVMT, ENPC, 6-8 avenue Blaise Pascal, Cit\'e Descartes
Champs-sur-Marne, 77455 Marne-la-Vall\'ee cedex 2, France.}
\email{frederic.meunier@enpc.fr}

\maketitle

 %\MSC
 %11E72, 16K20.
 %\EndMSC

 \smallskip

\begin{abstract}
The paper deals with $\sigma$-games on grid graphs (in dimension $2$
and more) and conditions under which any completely symmetric
configuration of lit vertices can be reached -- in particular the
completely lit configuration -- when starting with the all-unlit
configuration. The answer is complete in dimension 2. In dimension
$\geq 3$, the answer is complete for the $\sigma^+$-game, and for
the $\sigma^-$-game if at least one of the sizes is even. The case
$\sigma^-$, dimension $\geq 3$ and all sizes odd remains open.
\end{abstract}
%%---------------------Here the prologue ends---------------------------------

%%--------------------Here the manuscript starts------------------------------

%===================================================
%                     Introduction
%===================================================
\section*{Introduction}

A nice combinatorial game is the following. Suppose you have a graph
whose vertices can be {\em lit} or {\em unlit} (equivalently {\em
on} or {\em off}). When you push on a vertex, its state as well as
the state of its neighbors change. This kind of game is called a
{\em $\sigma^+$-game}.

You start with the all-off configuration. Can you find a sequence of
pushes such that you get a all-on configuration? The rather
unexpected answer is that it is always possible to find such a
sequence. Indeed Sutner proved \cite{Su87}

\begin{thm}[Sutner's theorem]\label{thm:sutner}
The all-on configuration can always be achieved starting from the
all-off configuration for a $\sigma^+$ game on any graph $G=(V,E)$.
\end{thm}

It is possible to define a similar game, the {\em $\sigma^-$-game},
for which pushing on a vertex changes the state of all its neighbors
but not its own state. In this case, things become harder since it
is not always possible to find a sequence achieving the all-on
configuration when starting from the all-off configuration. Simple
examples are provided by complete graphs with an odd number of
vertices, paths of odd length, etc.

\medskip

\begin{center}
\begin{figure}
\includegraphics[width=12cm]{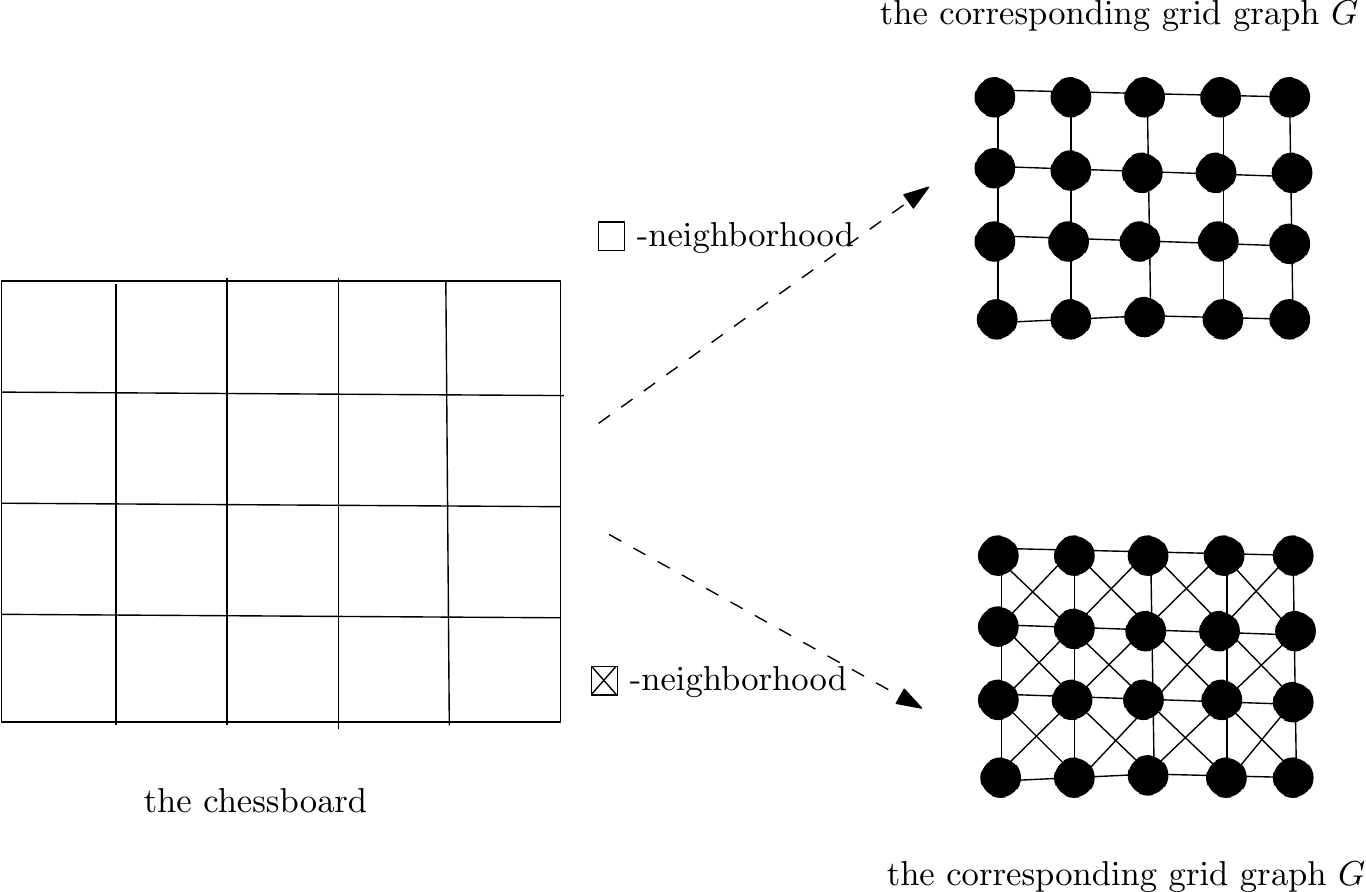}
\caption{The two usual kinds of neighborhood used for playing a
$\sigma$-game on a chessboard -- they lead to two distinct grid
graphs.} \label{fig:twokinds}
\end{figure}
\end{center}

$\sigma$-games have been intensively studied, and it not possible to
give here the whole list of references on this topic (see the
article \cite{GoWaWu09} for an extensive bibliography). Here we
focus on the case when the graph is a grid graph. Note that
$\sigma$-games on grid graphs have already been studied
(\cite{GoKl05} or \cite{BaRa96}, among many others) but for other
questions (for instance, the number of distinct configurations that
can be reached from a given one). Usually, two kinds of neighborhood
are considered for the grid graph: if the grid graph is seen as a
chessboard (the squares being the vertices), two squares sharing a
common edge are neighbors; depending whether two squares in contact
by their corners are or are not declared to be neighbors, we get one
or the other kind of neighborhood. The first kind of neighborhood is
denoted by $\square$ and the second one by $\boxtimes$. See Figure
\ref{fig:twokinds}. We will consider these kinds of neighborhood,
but also many others.

\medskip

In 2002, the French magazine ``Pour la Science'' published an
article written by Jean-Paul Delahaye and dealing with the
$\sigma^-$-game on grid graphs \cite{De02} (see also an updated
version of this article in the book \cite{De04}). The game was
defined on a chessboard and the neighbors of a square were the
adjacent squares, having a corner in common being enough to be
neighbor. Hence a square could have 8, 5 or 3 neighbors, depending
whether the square was or was not on the border or in the corner of
the chessboard (except if one dimension of the chessboard is 1, in
which case the number of neighbors is 2 or 1). We are here precisely
in the case of the $\boxtimes$-neighborhood.

In this article, a conjecture of a reader -- Nicolas Vaillant -- was
proposed. Recall that the 2-valuation of a number $n$ is the largest
$j$ such that $2^j$ is a divisor of $n$.

\medskip

\noindent{\bf Conjecture (Vaillant's conjecture):} The all-on
configuration cannot be achieved starting from the all-off
configuration for a $\sigma^-$-game played on a $m \times n$ chessboard if
and only if $n$ and $m$ are odd and such that $m+1$ and $n+1$ have
the same 2-valuation.

\medskip

In the present paper, we prove a general theorem (Theorem
\ref{principal}, Section \ref{sec:dim2}) that gives a necessary and
sufficient condition for a $\sigma$-game played on a $m \times n$
chessboard to be such that any doubly symmetric configuration can be
achieved. The approach will be purely algebraic. Vaillant's
conjecture is a consequence of this theorem.

As an other application, we obtain

\medskip

{\em Any doubly symmetric configuration can be achieved starting
from the all-off configuration for a $\sigma^+$-game played on a
chessboard for both the $\square$- and the
$\boxtimes$-neighborhoods.}

\medskip

A doubly symmetric configuration is a configuration that is
invariant by the symmetries with respect to the two medians of the
sides of the chessboard. Let us be more precise.

\begin{defi}
A configuration $Y=(y_{i,j})$ on a grid graph $m\times n$ is said to
be {\em doubly symmetric} if
$$y_{i,j}=y_{m+1-i,j}=y_{i,n+1-j}=y_{m+1-i,n+1-j}\quad\mbox{ for all
$i,j$.}$$
\end{defi}
The statements above (Vaillant's conjecture and the one concerning
doubly symmetric configurations for $\sigma^+$-games) are not only
true for the two usual kinds of neighborhood ($\square$ and
$\boxtimes$), but also for many others.

\medskip

The paper also deals with the case of $n_1\times
n_2\times\ldots\times n_d$ grids, where $d\geq 3$. We then speak about
completely symmetric configurations.

\medskip

\begin{defi}
A configuration $Y=((y_{j_1,\ldots,j_d}))$ on a $d$-dimensional grid
graph is said to be {\em completely symmetric} if
$$y_{j_1,\ldots,j_{i-1},j_i,j_{i+1},\ldots,j_d}=y_{j_1,\ldots,j_{i-1},n_i+1-j_i,j_{i+1},\ldots,j_d}\quad\mbox{
for all $i,j_1,j_2,\ldots,j_d$.}$$
\end{defi}

We will then prove the following result in Section \ref{sec:ddim}
(in a slightly more general form, Theorem \ref{thm:ddim}), but with
a different approach than that of Section \ref{sec:dim2}:

\medskip

{\em Any completely symmetric configuration can be achieved starting
from the all-off configuration for a $\sigma^+$-game played on a
$d$-dimensional grid for both the $\square$- and the
$\boxtimes$-neighborhoods..}

\medskip

It is also proved for many other kinds of neighborhoods. Note that
the question whether there is a simple condition for the existence
of a completely symmetric configuration when all dimensions are odd
for the $\sigma^-$-game remains unsettled (when $d\geq 3$, of
course...). Maybe this is due to the lack of an algebraic approach
for this case.

\section{Basic notions and notation}

Throughout this paper, we shall denote by $k$ a field. It will be of
characteristic $2$ starting from subsection \ref{Cheby}. We denote
by $\overline{\F}_2$ an algebraic closure of $\F_2$. For $n \geq 1$,
denote by $J_n$ the $n \times n$ matrix (with coefficients in $k$)
$$\begin{bmatrix}
  0      & 1 & 0 & \cdots & 0      \\
  1      & \ddots  & \ddots & 0 & \vdots \\
 0 & \ddots & \ddots & \ddots & 0\\
  \vdots & 0 & \ddots & \ddots & 1 \\
0 & \cdots & 0 & 1 & 0\\
\end{bmatrix}$$
 with $1$'s directly above and under the diagonal, and $0$'s everywhere else.

A \textit{game} $G$ on the $n \times m$ grid (the squares of which
can be lit or unlit) is given by the following. To each vertex $v$
of the grid, we associate a set of vertices whose state change if we
push on the vertex $v$. Equivalently, one may give a $nm \times nm$
matrix $M$ with coefficients in $\F_2$ (the field with two
elements), defined by the following property: let $1\leq i,k \leq n$
and $1 \leq j,l \leq m$ be integers. Then the coefficient of the
$(i,j)$-th column and the $(k,l)$-th line of $M$ is $0$ if pressing
on the vertex $(i,j)$ does not change the state of the vertex
$(k,l)$, and $1$ otherwise. We call $M$ the {\em generalized
adjacency matrix} of the game $G$. We will often assume that $M$
commutes with $J_n$ and $J_m$. The matrix $M$ can then be written as
a sum of scalar multiples of $J_n^r\otimes J_m^s$, where $r$ and $s$
are in $\mathbb{N}$ (cf. lemma \ref{commute}). For instance a $\sigma^-$-game played on a
$n\times m$ chessboard with a $\boxtimes$-neighborhood has a
generalized adjacency matrix $M=J_n\otimes J_m+J_n\otimes
I_m+I_n\otimes J_m$ (which of course commutes with $J_n$ and $J_m$).
Similarly, a $\sigma^+$-game played on a $n\times m$ chessboard with a
$\square$-neighborhood has a generalized adjacency matrix
$M=I_n\otimes I_m+J_n\otimes I_m+I_n\otimes J_m$ (which of course
commutes with $J_n$ and $J_m$).

All games will be assumed to be symmetric, i.e. that pressing on the
vertex $(i,j)$ changes the state of the vertex $(i',j')$ if and only
if pressing on the vertex $(i',j')$ changes the state of the vertex
$(i,j)$.

Using the generalized adjacency matrix $M$, to say that a
configuration can be achieved by the game $G$ starting from the
all-off configuration is equivalent to say that this configuration
(or more precisely the column vector of size $nm$ associated to it)
is in the image of $M$. We shall repeatedly use this elementary
remark without further mention.

We extend all these notions also for grids of dimensions $\geq 3$.

\medskip

Let us also recall how to prove Sutner's theorem (Theorem
\ref{thm:sutner}), which is valid for any kind of graph. The following lemma is a classical and elementary
result.

\begin{lem}\label{lem:binaryfarkas}
Let $U$ be a finite dimensional linear space over $k$, endowed with
a symmetric non-degenerate bilinear form, and let $\phi$ be a
self-adjoint endomorphism $U\rightarrow U$. We have
$$\mbox{\textup{Im} }\phi=\left(\mbox{\textup{Ker}
}\phi\right)^{\perp}.$$
\end{lem}

Theorem \ref{thm:sutner} is a straightforward consequence of Lemma
\ref{lem:binaryfarkas}, applied to $k=\F_2$, $U=\F_2 ^V$ and to
$\phi$ being the adjacency matrix of $G$ plus the identity matrix.
Indeed, it is then enough to prove that if we push on a subset $S$
of vertices that keeps the configuration in the all-off state, then
$S$ has cardinality even. But this is obvious since each vertex of
$G[S]$ must be of odd degree (otherwise some vertices of $S$ would
be on) and since the number of odd degree vertices in any graph is
always even.

We can reformulate this last sentence as a lemma, which will be
useful in the proof of Theorem \ref{thm:ddim}, in the particular
case of grid graphs. The matrix $M$ is the `generalized adjacency
matrix' of the game, defined in the beginning of this section.

\begin{lem}\label{lem:linsutn}
In the case of a $\sigma^+$-game, the number of nonzero entries of
any element of $\mbox{\textup{Ker }} M$ is even.
\end{lem}

\section{Doubly symmetric configurations on chessboards (or $2$-dimensional grids)}

\label{sec:dim2}

\subsection{Some algebra}

In this section, we introduce the technical material needed in the proof of the main theorem (Theorem \ref{principal}, Section \ref{sec:dim2}).\\
Let us begin with a key lemma.

\medskip

Let $A$ be a factorial ring, and $p \in A$ a prime element. For any
nonzero $x \in A$, we denote by $v_p(x)$ the highest power of $p$
dividing $x$.

\begin{lem}\label{lem:technic}
Let $p,q,r$ and $s$ be nonnegative integers. Consider the (local)
$k$-algebras $A=k[X]/X^p$ and $B=k[Y]/Y^q$. Denote by $x$ (resp. by
$y$) the class of $X$ (resp. of $Y$) in $A$ (resp. in $B$). In $A
\otimes_k B$, we still denote by $x$ the element $x \otimes 1$, and
similarly for $y$. Let $u= x-y\in A \otimes_k B$. Then the element
$x^r  y^s$ is divisible by $u$ if and only if $r+s \geq
\mbox{inf}\{p,q \}$. What is more, the same statement holds if we
replace $u$ by $cx+dy+\mbox{terms of order at least 2}$, where $c,d$
are nonzero elements of $k$.
\end{lem}

\begin{preuve}
Assume that $p \leq q$. Suppose that $r$ is positive. From the
relation $(x-y)(x^r  y^{s-1})=  x^{r+1} y^{s-1}-x^r y^s$, we deduce
that $x^r y^s$ is divisible by $x-y$ if and only if $x^{r+1}
y^{s-1}$ is so. Thus, we are reduced to the case where $s=0$. If $r
\geq \mbox{inf}
 \{p,q \}=p$,
then $x^r=0$, so one implication of the statement is obvious.
Conversely, assume that $x^r$ is divisible by $u$. Write
$x^r=(x-y)v$, for some $v$ in $A \otimes_k B$. We have a morphism
$$f: A \otimes_k B \longrightarrow A,$$ $$x \mapsto x,$$ $$y \mapsto
x.$$ Applying $f$ to the previous equality, we get $x^r=0$, hence
$r\geq p$, qed. For the last assertion, assume that $u=cx+dy+\lambda
x+\mu y$, where $\lambda, \mu$ lie in the maximal ideal $M$ of $A
\otimes_k B$. Clearly, we may assume that $c=d=1$. Put
$x'=x(1+\lambda)$ and $y'=-y(1+\mu)$. Then $u=x'-y'$ and $k[x',y']$
equals $A \otimes_k B$. Indeed, the obvious map $$f: A \otimes_k B
\longrightarrow k[x',y'],$$ $$x \mapsto x',$$ $$y \mapsto y'$$ is
injective, hence an isomorphism by dimension reasons. Note that
injectivity can be seen the following way: if $a$ is an element of
$M^n$ (where $n$ is a positive integer, and $M$ denotes as before
the maximal ideal of $A\otimes_k B$), we have $f(a)=a$ modulo
$M^{n+1}$. The second result of the lemma now follows by an
application of the preceding statement to $x'$ and $y'$. Indeed,
$x^r  y^s$ is divisible by $u$ if and only if  $x'^r  y'^s$ is so.

\end{preuve}\\

\noindent We can now state and prove the main proposition of this subsection.

\begin{prop} \label{linalg}
Assume that $k$ is algebraically closed. Let $P,Q,R,S$ be four
polynomials (in $k[X]$). Consider the $k$-algebras $A=k[X]/P$ and
$B=k[Y]/Q$. Let $U$ be an element of $k[X,Y]$. Put $u=U(x,y) \in A
\otimes_k B$. Assume the following: for every $\alpha, \beta \in k$
such that $P(\alpha)=Q(\beta)=0$ and $U(\alpha,\beta)=0$, we have
that $\frac {\partial U} {\partial X}(\alpha,\beta) \neq 0$ and that
$\frac {\partial U} {\partial Y}(\alpha,\beta) \neq 0$. Then $u$
divides $R(x)S(y)$ if and only if the following holds: for
\textit{every} $\alpha$, $\beta$ as above, denote by $p$ (resp. $q$,
$r$, $s$) the multiplicity of $\alpha$ (resp. $\beta$, $\alpha$,
$\beta$) as a root of $P$ (resp. $Q$, $R$, $S$). Then $r+s \geq
\mbox{inf}\{p,q \}$.

\end{prop}
\begin{preuve}
Write $P=(X-\alpha_1)^{m_1}...(X-\alpha_d)^{m_d}$. The Chinese
Remainder Theorem ensures that the natural morphism
$$A \longrightarrow k[X]/(X-\alpha_1)^{m_1} \times ... \times k[X]/(X-\alpha_d)^{m_d}$$
is an isomorphism. Using the similar isomorphism for $B$, we are
immediately reduced to the case where $P=(X-\alpha)^p$ and
$Q=(Y-\beta)^q$. If $U(\alpha,\beta) \neq 0$, then $u$ is invertible
in $A \otimes_k B$, hence the proposition is true in this case. If
$U(\alpha,\beta)=0 $, then replacing $P$ by $P(X+ \alpha)$ and $Q$
by $Q(Y+\beta)$, we may assume that $\alpha=\beta=0$, i.e. that
$P=X^p$ and $Q=Y^q$. We may also assume that $R$ and $S$ are powers
of $X$ (indeed, if $T$ is a polynomial such that $T(0) \neq 0$, then
$T(x)$ (resp. $T(y)$) is invertible in $A$ (resp. in $B$)). The
content of the proposition then boils down to that of lemma
\ref{lem:technic}, since the hypothesis about partial derivatives
ensures that $u$ is of the form $cx+dy+ \mbox{higher order terms}$.
\end{preuve}

\begin{lem} \label{drond}
Let $P,Q \in k[X]$ be two polynomials. Put $A=k[X]/P(X)$ and
$B=k[Y]/Q(Y)$. Denote by $x$ (resp. by $y$) the class of $X$ (resp.
of $Y$) in $A$ (resp. in $B$). Let $u$ be an element of $A\otimes_k
B$. Let $U \in k[X,Y]$ be a polynomial such that $u=U(x,y)$. Assume
that $\alpha \in k$ is a root of multiplicity $\geq 2$ of $P$, and
let $\beta \in k$ be any root of $Q$. Then the partial derivative
$\frac {\partial U} {\partial X} (\alpha, \beta)$ is independent of
the choice of $U$.
\end{lem}

\begin{preuve}
Indeed, any other $U'$ satisfying $u=U'(x,y)$ is of the form
$U'=U+R(X,Y)P(X)+S(X,Y)Q(Y)$, and the hypothesis about $\alpha$ implies
that   $\frac {\partial U} {\partial X} (\alpha, \beta) =\frac
{\partial U'} {\partial X} (\alpha, \beta)$.
\end{preuve}

\begin{defi} \label{partindep}
Under the hypothesis of lemma \ref{drond}, we shall denote  $\frac
{\partial U} {\partial X} (\alpha, \beta)$, which is independent of
the choice of $U$, by $\frac {\partial u} {\partial x} (\alpha,
\beta)$.
\end{defi}

\subsection{Preliminaries on Chebychev polynomials} \label{Cheby}

Chebychev polynomials modulo 2 are classical tools in the study of
$\sigma$-games on grid graphs in dimension 2 (see \cite{Su00} for
instance, or \cite{GoKl05}, where they are called Fibonacci
polynomials). We recall in this subsection their definition and some
of their properties.

\subsubsection{Classical Chebychev polynomials}

\noindent The usual Chebychev polynomials are elements of $\mathbb Z [X]$ defined as follows.\\
\noindent Set $P_0=2$ and $P_1=X$. Then, define $P_n$ inductively by the formula
$$P_{n+1}=XP_n+P_{n-1}. $$
This formula will be called the \textit{Chebychev relation.}\\
\noindent The $P_n$'s satisfy the following well-known properties, valid for all nonnegative integers $n$ and $m$:\\

\noindent i) $P_n(X+X^{-1})=X^n+X^{-n},$\\
\noindent ii) $P_n P_m =P_{n+m}+P_{|n-m|}$.\\

\noindent Property i) in fact characterizes the Chebychev polynomials, and ii) is an easy consequence of i).

\subsubsection{Chebychev polynomials modulo $2$}

From now on, $k$ will be assumed to have characteristic $2$. For a
nonnegative integer $n$, define $Q_n$ to be the class of $\frac
{P_{n+1}} X$ in $k[X]$. It is readily seen that the $Q_n$'s are
indeed polynomials since all $P_n$'s are divisible by $X$ modulo
$2$. Note that the $Q_n$'s also satisfy the Chebychev relation. We
shall now study some elementary divisibility properties of these
polynomials. First of all, an easy induction shows that $Q_n$ is
divisible by $X$ if and only if $n$ is odd. From point ii) of the
preceding section, we have that $Q_{2n-1}=XQ_{n-1}^2$. This implies
a formula useful in the sequel. Take an odd positive integer n.
Write $n+1=2^j m$, with $m$ odd. The preceding formula, applied
several times, then yields $Q_n=X^{2^j-1}Q_{m-1}^{2^j}$. Since $m-1$
is even, we have that $Q_{m-1}$ is not divisible by $X$, hence the
relation:
$$v_X(Q_n)=2^j-1.$$
We also get the following. Let $n$ be odd, and $R \neq X$ be a
(monic) prime polynomial dividing $Q_n$. From the relation
$Q_{n}=XQ_{\frac {n-1} 2}^2$, we infer that

$$v_R(Q_n)=2v_R(Q_{\frac {n-1} 2})$$
and
$$v_X(Q_n)=1+2v_X(Q_{\frac {n-1} 2}).$$

\noindent Those relations are basically the only facts we shall need about Chebychev polynomials.

\subsection{Statement and proof of the main theorem}

It is an elementary exercise to check that the characteristic
polynomial of $J_n$ is $Q_n$. Let $e_i$ ($i=0... n-1$) denote the
$i$'th basis vector of $k^n$. Consider the linear map

$$\Phi_n: k[X]/Q_n \longrightarrow k^n, $$
$$ X^i \mapsto J_n^i(e_1). $$

This map is well-defined (Cayley-Hamilton). It is readily checked
that it is surjective, hence an isomorphism (this amounts to saying that the characteristic and minimal polynomials of $J_n$ coincide). In the sequel, we will
\textit{identify} $k^n$ with $k[X]/Q_n$ using $\Phi_n$. We shall
denote by $x_n$ the class of $X$ in $k[X]/Q_n$. One sees that, under
the isomorphism given by $\Phi_n$, $e_i$ corresponds to $Q_i(x_n)$.
Furthermore, the action of $J_n$ on $k[X]/Q_n$ is simply given by
multiplication by $x_n$.

\begin{lem} \label{commute}
Let $f$ be an endomorphism of the $k$-vector space $k[X]/Q_n
\otimes_k k[Y]/Q_m$ commuting with multiplication by $x_n=\overline
X \otimes 1$ and $y_m=1 \otimes \overline Y$. Then $f$ is given by
multiplication by $f(1)$.

\end{lem}
\begin{preuve}
Easy verification.
\end{preuve}

\begin{defi}(central configuration)
 The central configuration of the $n \times m$ grid is defined the following way: put $c_n:=Q_{\frac {n-1} 2}(x_n)$ if $n$ is odd, $c_n:=Q_{\frac {n} 2}(x_n) + Q_{\frac {n} 2 -1}(x_n)$ if $n$ is even, and define $d_m$ similarly with respect to $y_m$. Then the central configuration is $c=c_n d_m$. It consists in the central square if $n$ and $m$ are both odd, in the $2$ central squares if exactly one of the two integers $n$ and $m$ is odd, and in the $4$ central squares if $n$ and $m$ are both even.
\end{defi}

\begin{lem} \label{reduc}
 Every doubly symmetric configuration in $k[X]/Q_n \otimes_k k[X]/Q_m$ is divisible by the central one. Moreover, if $n$ and $m$ are both odd, then the central configuration is divisible by the all-on configuration.
\end{lem}

\begin{preuve}
Let us prove the first assertion. It suffices to show that $e_{n-1-i} + e_{i}=Q_{n-1-i}(x_n) + Q_{i}(x_n)$ is divisible by $c_n$ in $k[X]/Q_n(X)$ for any $n\geq 1$ and any integer $i$ satisfying $0 \leq i \leq [n/2]$. This is an easy descending induction on $i$, using the relation $Q_{k+1}=XQ_k+Q_{k-1}$. Let us now handle the second assertion. It suffices to prove it for the 1-dimensional case. The all-on configuration is then $\sum_{i=0 \ldots n-1} Q_i(x_n) = Q_{\frac {n-1} 2 } (x_n)( Q_{\frac {n-1} 2 } (x_n)+  Q_{\frac {n-3} 2 } (x_n))=c_n(c_n +  Q_{\frac {n-3} 2 }(x_n))$ by a straightforward computation. It is enough to show that $ Q_{\frac {n-1} 2 }(x_n) + Q_{\frac {n-3} 2 }(x_n) $ is invertible in $k[X]/Q_n(X)$, i.e. that $ Q_{\frac {n-1} 2 }(X)+ Q_{\frac {n-3} 2 }(X)$ and $Q_n(X)=X Q_{\frac {n-1} 2 } ^2(X)$ are coprime. Let $R$ be a monic prime polynomial dividing these two polynomials. Certainly $R$ is not $X$ since the constant term of $ Q_{\frac {n-1} 2 }(X)+ Q_{\frac {n-3} 2 }(X)$ is $1$. But then $R$ divides both $ Q_{\frac {n-1} 2 }(X)$ and $Q_{\frac {n-3} 2 }(X)$, which are coprime.
\end{preuve}

\begin{thm} \label{principal}
Let $n$ and $m$ be integers. Let $G$ be a game on the $n \times m$
grid (identified with $\F_2 [X]/Q_n \otimes_{\F_2} \F_2 [Y]/Q_m$)
which commutes with the elementary games $J_n$ and $J_m$. Let $f$ be
the endomorphism of  $\F_2[X]/Q_n \otimes_{\F_2} \F_2[Y]/Q_m$ given
by $G$- it corresponds to the generalized adjacency matrix $M$ of
$G$. By Lemma \ref{commute}, $f$ is then given by multiplication by
$u:=f(1)$. For any $\alpha \in \overline{\F}_2$ (resp. $\beta \in
\overline{\F}_2$), which is a root of $Q_n$ (resp.
$Q_m$), such that $u(\alpha, \beta)=0$, we assume the following. \\
If $\alpha$ (resp. $\beta$) is a root of multiplicity $\geq 2$ of
$Q_n$ (resp. $Q_m$), then $\frac
{\partial u} {\partial x} (\alpha, \beta) \neq 0$ (resp. $\frac
{\partial u} {\partial y} (\alpha, \beta) \neq 0$) (these quantities
are well-defined thanks to Lemma \ref{drond}). \\ Then any doubly
symmetric configuration can be achieved starting from the all-off
configuration for the game $G$ if and only if the three following
conditions do not simultaneously hold:  $n$ and $m$ are both odd,
$u(0,0)=0$ and $v_2(n+1) = v_2(m+1)$.

\end{thm}
\begin{rem}
In the case where $n$ and $m$ are both odd, the fact that any doubly
symmetric configuration can be achieved is equivalent to the fact
that the all-on configuration can be achieved; this is the content
of Lemma \ref{reduc}.
\end{rem}

\begin{rem}
Let us be more precise concerning how the theorem implies Vaillant's
conjecture, stated in the beginning of the paper. This is the
$\boxtimes$-case
$$u=y_m+x_n+y_m x_n.$$  The condition about partial derivatives is
here obvious: indeed, put $U=X+Y+XY$. We have $\frac {\partial U}
{\partial X}=1+Y$, so that the condition may fail only for
$\beta=1$. But $U(X,1)=1$ has no root in $k$.

For the $\square$-case $$u:=y_m+x_n,$$ we see that the conclusion is
identical.

In the case of a $\sigma^+$-game with the usual neighborhoods
defined in the introduction (the $\square$ and the $\boxtimes$
neighborhoods), i.e. when $u:=y_m+x_n+1$ or $u=y_m+x_n+y_m x_n+1$,
we have $u(0,0)=1$. Hence, any doubly symmetric configuration can
always be achieved.
\end{rem}

\medskip

\begin{preuve}[Proof of Theorem \ref{principal}]
Put $k=\F_2$. Let $A=k[X]/Q_n$ and $B=k[Y]/Q_m$. By lemma
\ref{reduc}, it is enough to show that the central configuration can
be obtained. Over fields, the formation of the image of a linear map
commutes with scalar extension. Thus, there is no harm in replacing
$k$ by an algebraic closure of $k$. Let $U \in k[X,Y]$ be such that
$U(x,y)=u$.  Assume the hypothesis about partial derivatives holds.
If $\alpha$ is a simple root of $Q_n$, then $Q_n'(\alpha) \neq 0$.
Because $k$ is infinite, we can then replace $U$ by $U + \lambda
Q_n(X)$, for a suitable $\lambda \in k$, in such a way that the
partial derivatives of $U$ with respect to $X$ are nonzero when
evaluated at $(\alpha, \beta)$, where $\alpha$ runs through the
simple roots of $Q_n$ and $\beta$ through the roots of $Q_m$,
submitted to the condition $u(\alpha,\beta)=0$. We can then do the
same for partial derivatives with respect to $Y$. By doing so, we
get a $U \in k[X,Y]$ such that $U(x,y)=u$ and   $\frac {\partial U}
{\partial X}(\alpha,\beta) \neq 0$,  $\frac {\partial U} {\partial
Y}(\alpha,\beta) \neq 0$ for any $\alpha \in k$ (resp. $\beta$)
which is a root of $Q_n$ (resp. $Q_m$) and such that $u(\alpha,
\beta)=0$. We now want to apply Proposition \ref{linalg}. Let
$\alpha$, $\beta \in k$ be as before. Put $P:=Q_n$, $Q:=Q_m$,
$p:=v_{X-\alpha}(P)$, $q:=v_{Y-\beta}(Q)$. Put $R:=Q_{\frac {n}
2}(X) + Q_{\frac {n} 2 -1}(X)$ if $n$ is even, $R:=Q_{\frac {n-1}
2}(X)$ if $n$ is odd. Similarly, put $S:=Q_{\frac {m} 2}(Y) +
Q_{\frac {m} 2 -1}(Y)$ if $m$ is even, $S:=Q_{\frac {m-1} 2}(Y)$ if
$m$ is odd. Put $r:=v_{X-\alpha}(R)$ and $s:=v_{Y-\beta}(S)$. If $n$
is even, we compute: $$XR^2=XQ_{\frac {n} 2}(X)^2 + XQ_{\frac {n} 2
-1}(X)^2=$$
$$Q_{n+1}(X)+Q_{n-1}(X)=XQ_n(X).$$ Hence $p=2r$. If $n$ is odd, we
then have $1+2r=p$ if $\alpha=0$ and $2r=p$ otherwise, as proved in
section \ref{Cheby}. Similarly, we get relations between $s$ and
$q$. To finish the proof, we have to show that the relation $r+s
\geq \mbox{inf}\{p,q \}$ (for every $\alpha$ and $\beta$) is
equivalent to the fact that the three conditions of the theorem do
not simultaneously hold. Assume that $n$ is even. If $s \geq r$,
then $r+s \geq 2r=\mbox{inf}\{p,q \}$. If $s < r$, then $r+s \geq
2s+1 \geq q=\mbox{inf}\{p,q \}$. Hence the relation is valid in this
case. In the same way, it is valid if $m$ is even. Assume now that
$n$ and $m$ are both odd. If $\alpha \neq 0$, then $2r=p$, and we
conclude as before that the relation holds. In the same way, it
holds if $\beta \neq 0$. Assume now that $\alpha=\beta=0$ (hence
that $u(0,0)=0$). Then $1+2r=p$ and $1+2s=q$, hence the relation
$r+s \geq \mbox{inf}\{p,q \}$ holds if and only if $p$ and $q$ are
distinct, which in view of Section \ref{Cheby} amounts to saying
that $v_2(n+1) \neq v_2(m+1)$.

\end{preuve}

\section{Completely symmetric configurations on $d$-dimensional grids}

\label{sec:ddim}

We extend in this section some of the previous results.

\begin{thm} \label{thm:ddim}
Let $G$ be a game on the $n_1\times\ldots\times n_d$ grid that
commutes with the elementary games $J_{n_i}$ for $i=1,\ldots,d$. If
$G$ is a $\sigma^+$-game, then any completely symmetric
configuration can be achieved starting from the all-off
configuration.
\end{thm}

Note that the dimension 2 case is also covered, giving an
alternative proof of some statements already proved in the previous
section for $\sigma^+$-games.

The following lemma plays a crucial role in the proof. In a sense,
it proves the theorem above for the 1-dimensional case. Denote by
$S$ the map
$$\begin{array}{cccl}\mathbb{F}_2^n & \longrightarrow & \mathbb{F}_2^{\lceil n/2\rceil} & \\
(y_1,\ldots,y_n) & \longmapsto & \begin{array}{cl}(y_1+y_n,y_2+y_{n-1},\ldots,y_{n/2-1}+y_{n/2}) & \text{ if $n$ is even,} \\ (y_1+y_n,y_2+y_{n-1},\ldots,y_{\lceil n/2\rceil -1}+y_{\lceil n/2\rceil +1},y_{\lceil n/2\rceil}) & \text{ if $n$ is odd,} \end{array} & \end{array}$$
and by $c$ the map $$\begin{array}{ccc}\mathbb{F}_2^n & \rightarrow & \mathbb{F}_2\\
(y_1,\ldots,y_n) & \mapsto & \sum_{i=1}^ny_i \end{array}$$

\begin{lem}\label{lem:symlin}
Let $F$ be a linear subspace of $\mathbb{F}_2^n$ such that
$J_nF\subseteq F$ and $F\subseteq\mbox{\textup{Ker }} c$. Then $F\subseteq\mbox{\textup{Ker }}S$.\end{lem}

$F\subseteq\mbox{\textup{Ker }} c$ reads also {\em any element of
$F$ has an even number of nonzero entries}. Note the similarity with
the statement of Lemma \ref{lem:linsutn}. Indeed, we will apply
Lemma \ref{lem:symlin} to some linear spaces constructed from $M$.

\medskip

\begin{preuve}[Proof of Lemma \ref{lem:symlin}]
Suppose  that there is $\boldsymbol{x}$ in $F$ with a $i$ such that
$x_i\neq x_{n+1-i}$. Since $J_n^r\boldsymbol{x}\in F$ for any
positive integer $r$, we can assume that $i=1$.

We have $J_n\boldsymbol{x}\in F$, hence $c( J_n\boldsymbol{x})=0$ and
$$x_2+(x_1+x_3)+\ldots+(x_{n-2}+x_n)+x_{n-1}=0,$$ whence $x_1+x_n=0$, a contradiction.

Therefore, for all $\boldsymbol{x}\in F$ and all $i=1,\ldots,n$, we have $x_i=x_{n+1-i}$.

When $n$ is odd, $x_{\frac{n+1}2}=0$ is then a direct consequence of
the fact that any element of $F$ has an even number of nonzero entries.

\end{preuve}

We now restate Lemma \ref{lem:symlin} in a slightly more general
form, which will suit the proof scheme of Theorem \ref{thm:ddim}.

\begin{lem}\label{lem:symlind}
Let $V$ be any $\mathbb F_2$-vector space. Let $F$ be a linear
subspace of $V\otimes \mathbb{F}_2^n$ such that $(Id \otimes
J_n)(F)\subseteq F$ and $F\subseteq\mbox{\textup{Ker }} (Id \otimes
c)$. Then $(Id \otimes S)(F) \subseteq V\otimes\mathbb{F}_2^{\lceil
n/2\rceil} $ is zero.

\end{lem}

\begin{preuve}
Let  $\phi$ be any linear form on $V$. The lemma is a direct
consequence of Lemma \ref{lem:symlin} once we have noticed that
$(\phi \otimes Id)(F) \subseteq \F_2^n$ satisfies the conditions of
Lemma \ref{lem:symlin} for $F$, and that the following diagram
commutes: $$ \xymatrix {V\otimes \F_2^n  \ar[d]^{\phi \otimes Id}
\ar[r]^{Id \otimes S } & V\otimes \F_2^{\lceil n/2\rceil}
\ar[d]^{\phi \otimes Id} \\ \F_2^n \ar[r]^S & \F_2^{\lceil
n/2\rceil}.} $$

The vector space $(Id \otimes S)(F)$ is thus killed by all linear
forms on $V$, hence is zero.

\end{preuve}
\\
\begin{preuve}[Proof of Theorem \ref{thm:ddim}]
Let $M$ be the adjacency matrix of the game $G$. The hypothesis that
$M$ commutes with all $J_{n_i}$'s ensures that $M$ is a linear
combination of matrices of the form $J_{n_1}^{i_1} \otimes\ldots
\otimes J_{n_d}^{i_d}$- the proof of this fact is the same as that
of lemma \ref{commute}.
%It is thus self-adjoint for the standard
%nondegenerate symmetric bilinear form on $\F_2 ^ G = \F_2 ^ {n_1}
%\otimes \ldots \otimes \F_2 ^ {n_d}$.
Define $N:=\mbox{Ker }M$. Note
that $N$ is stable by all $J_{n_i}$. Denote by $S_i : \F_2^{n_i}
\longrightarrow \F_2^{\lceil n_i/2\rceil}$ (resp. $c_i : \F_2^{n_i}
\longrightarrow \F_2$) the map defined in the same way as $S$ (resp.
$c$), for $n=n_i$.

We have the following property.

\medskip

{\em Let $i\in\{0,1,\ldots, d\}$. Then $$S_1\otimes S_2\otimes \ldots \otimes S_i \otimes c_{i+1} \otimes c_{i+2}\otimes \ldots\otimes c_d (N)=\{0\}.$$}

Indeed, it is true for $i=0$, according to Lemma \ref{lem:linsutn}.

The other cases are true by induction. Suppose that the property is true for $i\geq 0$. Define
 $$F:=S_1\otimes S_2\otimes\ldots\otimes S_i\otimes Id \otimes c_{i+2}\otimes c_{i+3}\otimes \ldots\otimes c_d (N).$$
$F$ is a linear subspace of $V\otimes\mathbb{F}_2^{n_{i+1}}$ where
$V=\mathbb{F}_2^{\lceil
n_1\rceil}\otimes\ldots\otimes\mathbb{F}_2^{\lceil n_i\rceil}$. We
have $(\mathop{Id}_V \otimes J_{n_{i+1}})F\subseteq F$, and, by
induction, $(Id_V \otimes c_{i+1})(F)=\{0\}$. We apply Lemma
\ref{lem:symlind} and get that $(Id_V \otimes S_{i+1})(F)=\{0\}$,
which means exactly
$$S_1\otimes S_2\otimes\ldots\otimes S_{i+1}\otimes c_{i+2}\otimes
c_{i+3}\otimes \ldots\otimes c_d (N)=\{0\}.$$

\medskip

The statement of the theorem for the $\sigma^+$-game is a direct
consequence of the property above for $i=d$: apply Lemma
\ref{lem:binaryfarkas} to get that any completely symmetric
configuration is in the image of $M$.

\end{preuve}

The previous theorem has a nice corollary concerning
$\sigma^-$-games in any dimension.

\begin{coro}
Let $G$ be a $d$-dimensional $\sigma^-$-game on a
$n_1\times\ldots\times n_d$ grid, with $n_1$ even. Assume it can be
written as a finite sum
$\sum_{i_1,\ldots,i_d}\lambda_{i_1,\ldots,i_d}J_1^{i_1}\otimes\ldots\otimes
J_d^{i_d}$ ($\lambda_{i_1,\ldots,i_d} \in \F_2$), where
$\lambda_{1,0,\ldots,0}=1$ and where all other
$\lambda_{i_1,\ldots,i_d}$ equal $0$ except possibly when at least
one of the $i_j$, $j \geq 2$, equals $1$. Then every completely
symmetric configuration can be achieved starting from the all-off
configuration.
\end{coro}

\begin{preuve}
Define the game $M':=J_{n_1}^{-1}M$. The hypothesis of the corollary
ensures that $M'$ is the generalized adjacency matrix of a
$\sigma^+$-game. Theorem \ref{thm:ddim} applies to $M'$. The result
follows, for the space of completely symmetric configurations is
stable by $J_{n_1}$.
\end{preuve}

\bibliographystyle{amsplain}
\bibliography{sigma}
\end{document}